# Pole Shape Optimization Using a Tabu Search Scheme

P. J. Leonard and A. M. Connor

***Abstract:*** **The pole shape optimization of an electromagnet typical of an MRI type application is investigated. We compare the use of different parameterizations of the pole shape and compare the convergence of the optimizations using a discrete variable step length Tabu Search scheme.**

***Index Terms:*** **Magnetic resonance imaging, magnetostatics, optimization methods, pole shape, Tabu.**

## I. Introduction

MRI applications require a high field homogeneity within a region. If iron poles are used these can be shaped to increase the field uniformity. New powdered iron materials can be molded into complex contours [1]. Choosing a good parameterization of the shape is critical to the success of any optimization process. The designer faces a dilemma:

- Choose a detailed representation using many parameters and risk the optimization process running out of time before a good solution is found.
- Choose a coarse representation with only a few parameters and risk the optimal solution being outside the search space.

In the case of pole shape optimization experience has shown that a good "rule of thumb" is a profile with 2 steps [2]. In this paper we will investigate pole shape optimization using different parameterizations representing different numbers of steps in the pole.

## II. Tabu Search Scheme

The basic concept of Tabu search as described by Glover [3] is a metaheuristic superimposed on another heuristic. The overall approach is to adopt an aggressive strategy that forces the underlying heuristic to always make a move, even if the search is trapped in a local optima. A move is a transformation from one solution to a new solution in the allowable neighborhood.

Following a steepest descent/mildest ascent approach, a move may either result in a best possible improvement or a least possible deterioration of the objective function value. Without additional control, however, such a process can cause a locally optimal solution to be re-visited immediately after moving to a neighboring solution which then results in an unending cycle between the two solutions.

To prevent this cycling, Tabu search introduces the concept of an attribute based memory. By choosing suitable attributes of either a move or the solution resulting from a move it is possible to retain a representative record of the search trajectory through the solutions space. This representative record can be used to guide the search trajectory away from local optima in an attempt to explore different regions of the solution space.

### A. Implementation

The principles of Tabu search can be implemented in a number of ways depending on a combination of the underlying search algorithm, the attributes chosen to be used to represent the search vector through the solution space and the enhancements added to control the search. The method described in this paper is a very simple implementation where the attributes chosen are simply the values of the design parameters and the underlying search algorithm is based on a Hooke and Jeeves method [4]. Simple search intensification and diversification strategies have been designed to alter the focus of the search vector without adding computational expense (see Figure 1).

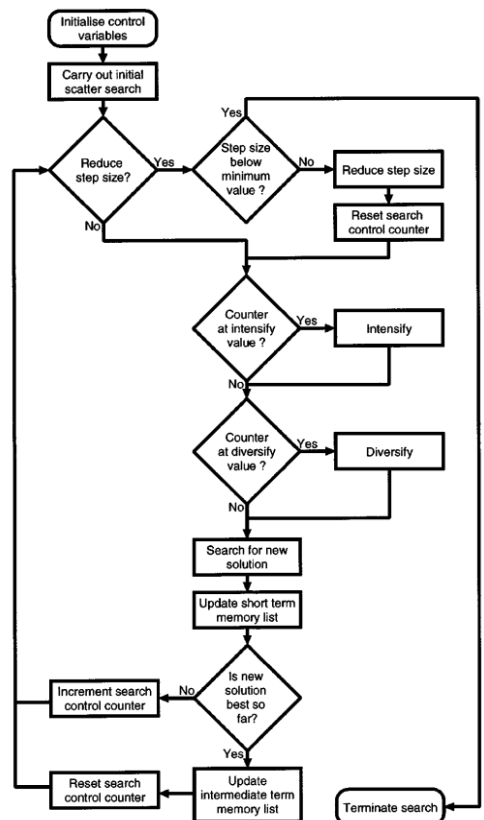

**Figure 1.** Tabu control algorithm.

## B. Short Term Memory

In general, the short term memory consists of a list of solution attributes that cannot be embodied in moves whilst they are still contained in the memory. The search is generally barred from making a move that essentially reverses a move contained in the list. This approach has the advantage that intensification and diversification of the search can be easily controlled by considering which moves have generally resulted in improvement and which moves are not regularly used. However, the storing of moves does not guarantee that the method will not return to the same real solution. As parametric design problems tend to deal with numerical design parameters rather than symbolic parameters the attributes used in the short term memory are simply the values of the design parameters themselves. This prevents the method returning to a previously visited solution by a convoluted route for as long as the solution is maintained in the list. The short term memory is continually refreshed as the search explores the solution space. As each move is made a new set of attributes are placed in the memory and the oldest set of attributes are removed.

## C. Search Intensification and Diversification

The short term memory enables the method to leave locally optimal solutions in the quest for the global optimum of a function. However, short term memory alone does not ensure that the search will be both efficient and effective. Search intensification and diversification techniques are often used first to focus the search in particular areas and then to expand the search to new areas of the solution space. This is normally achieved by the use of longer term memory cycles. Given the simple nature of the attributes used, the intermediate memory cycle is nothing more than a list of previously visited solutions that have resulted in an improved solution being found.

Without effective intensification and diversification strategies, the Tabu search method offers little over methods such as Dynamic Hill Climbing [5]. The work in this paper utilizes a very simple intensification strategy that proposes a new solution based on the average value for each design parameter calculated from the solutions contained in the intermediate memory list. During the initial stages of the search this often leads to a new good solution that is significantly different from those contained in the intermediate memory list as those solutions are still quite distinct. At later stages of the search this strategy tends to produce solutions that are very similar to those contained in the intermediate memory list as the method will have converged to a region in the solution space.

Previous implementations of the Tabu search method have used simple random refreshment of the search point to a new solution. The method used in this paper is working toward a more effective diversification strategy by introducing concepts of scatter search. Rather than use a single refreshment point, the strategy proposes a number of solutions and starts searching from the best solution found. At current, there is no control over the spacing of the solutions generated in the diversification though a potential refinement would be to introduce such control to guarantee that the solutions are suitably spaced throughout the solution space.

### D. Hill Climbing Algorithm

The underlying hill climbing algorithm used in this work is based upon the method developed by Hooke and Jeeves [4]. This method consists of two stages, the first of which carries out an initial exploration around a given base point. When a move to a new point (the exploration point) which improves the objective function is identified, the search is extended along the same vector by a factor in a pattern move. If this new solution has a better objective value than the exploration point, then this point is used as the new base point and the search is repeated. Otherwise, the search is repeated using the exploration point as the new base point. The algorithm used in this implementation of Tabu search differs in several ways from the standard Hooke and Jeeves algorithm.

In the Hooke and Jeeves algorithm each parameter is varied in turn and the first move that results in a better objective function value is selected. The implication of this is that not all potential moves are evaluated. In the Tabu search implementation all trial moves are investigated and the best move that is not tabu is chosen.

Another difference concerns how the step size is periodically reduced. In the Hooke and Jeeves search when a point is reached from which no improvement can be found then the step size is reduced by a factor of two. In the Tabu search implementation this is not practical as the Tabu search metaheuristic forces the search point out of local optima. The step size is therefore reduced when other conditions apply. A counter is maintained of the number of search moves that have elapsed since an improved solution was found. When this reaches a given value then the search carries out an intensification action. If an improved solution is found then the counter is reset. If no improvement is found then the search continues until the counter reaches a higher preset value at which point diversification is carried out. Again, the search continues and if no improvement is found before the number of moves reaches the next preset level then the step size is reduced.

At present there are no convergence criteria built into the search control to terminate the search before the step sizes are brought to their minimum sizes and the search ends naturally. This implies that a potentially large number of evaluations are used to produce only very small improvements in objective function value. A more efficient search could easily be produced by introducing convergence criteria.

We used the same control strategy that had been found to be successful for previous problems. If the dimension of the parameters is $N$ then the intensification count is max(4, $N$/2), diversification count is max(8, $N$) and the reduce step size count is max(3$N$/2, 12) though it is acknowledged that possible other values may have yielded more efficient optimization processes.

## III. RESULTS

An idealized pole shape problem was investigated (see Figure 2). The configuration is typical of a MRI application. A pair of rotationally symmetric poles are used to shape the field.

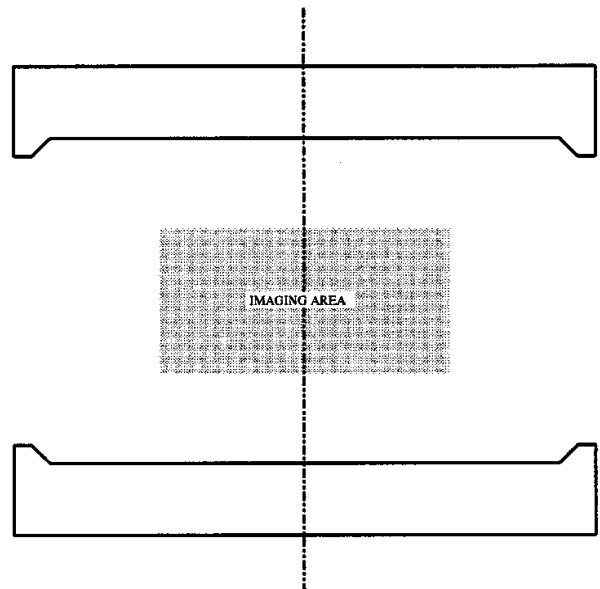

**Figure 2.** Idealized MRI application.

For simplicity the source of the field is not modeled, the poles are assumed to driven by auxiliary coils or permanent magnets, that set up a m.m.f between the pole faces.

This was modeled using two dimensional axisymmetric finite elements. Symmetry allows one half of the axisymmetric problem to be modeled. A magnetic scalar formulation was used with the pole set to a fixed potential. The finite element model was truncated with a natural boundary at a radius of 400 mm and a distance of 300 mm from the plane of symmetry.

*A. Detailed Parameterization* RAW

The pole face shape was parameterized with a 16 segment line (these corresponded to the edges of finite elements). The points were evenly spaced in the radial direction but the axial coordinate was allowed to vary within the range 0–40 mm.

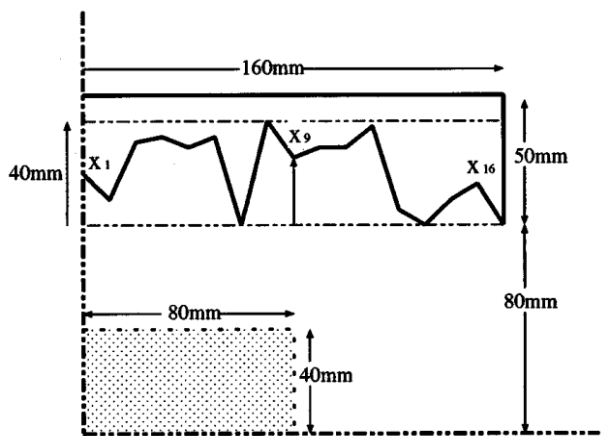

**Figure 3.** The axisymmetric finite element model.

The elements 50 mm above and 40 mm below this line where stretched and compressed to accommodate this geometry change, this avoids numerical discontinuities that can arise if remeshing is used. We looked at the field homogeneity in a cylinder of radius 80 mm and height of 80 mm. The $z$ directed field, $B = B_Z$ was sampled in the target region. The cost function was defined as cost = $(B_{MAX} - B_{MIN}) / (B_{MAX} + B_{MIN})$.

The problem space was made discrete by choosing a resolution of 0.00001 m. The initial step size was chosen to be 1/5 of the parameter range. It is known (from previous experience) that an initial flat pole shape is a good starting point for an optimization. However, to make life more interesting for the optimization scheme we denied it the privilege of this information. We started at random points in the search space.

*B. Reduced Parameter Representations*
We investigated reducing the number of parameters that represented the pole shape. Four cases were investigated, representing 1, 2, 3 and 4 ramps in the pole height. Figure 4, shows the 2 ramp version. The radial length parameters $a$ and $b$ were given a minimum step of 10 mm, (corresponding to node positions).

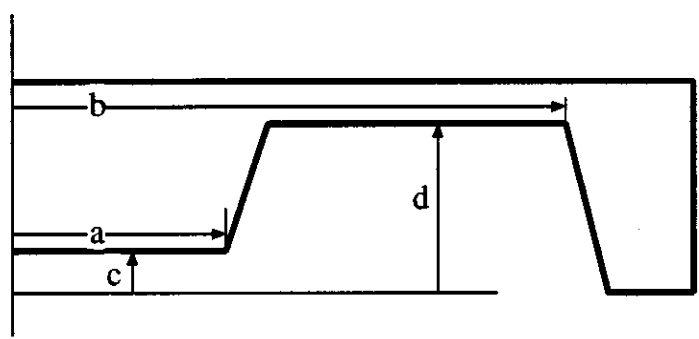

**Figure 4.** Four parameter study.

The parameter $b$ was allowed to take the range 10mm $<= b <=$ 160mm, whilst $a$ was constrained so $0 <= a < b$. The 2 heights $c$ and $d$ used the same resolution as the RAW parameter version, making this scheme a subspace of RAW. The 1, 3 and 4 step schemes were constructed in a likewise fashion.

*C. Results*

It is difficult to compare schemes when they involve an element of random choice. We choose to perform five optimizations for each parameterization. For each parameterization we noted the objective function after 20, 100 and 1000 iterations, we also note the final value. The results are displayed in Table I which shows the range of observed values for each parameterization.

**Table I.** Variation of the objective function versus iterations over 5 runs using different parameterizations

| Scheme | Iteration / objective function $\times 10^3$ | | | |
|---|---|---|---|---|
| | *20* | *100* | *1000* | *inf* |
| RAW | 3029 - 3395 | 1701 - 2534 | 172 - 1999 | 112 - 315 |
| 1 | 1006 - 3897 | 597 - 1118 | 520 - 595 | 520 - 595 |
| 2 | 610 - 2252 | 284 - 1136 | 126 - 250 | 112 - 149 |
| 3 | 1105 - 4474 | 213 - 2343 | 66 - 425 | 55 - 217 |
| 4 | 6065 - 13256 | 576 - 1419 | 183 - 429 | 98 - 422 |

Some typical convergence characteristics are shown in Figs. 5 and 6. Typical pole shapes found are shown in Figure 7.

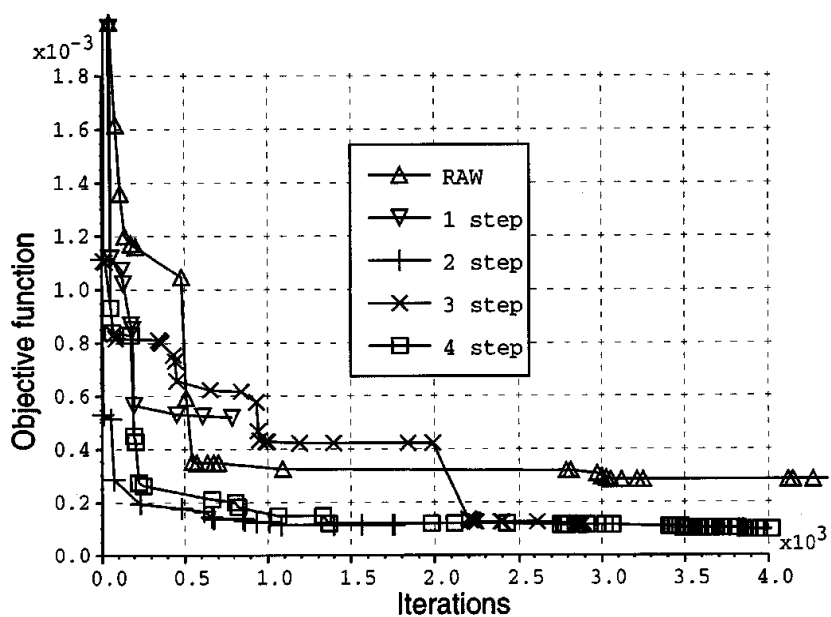

**Figure 5.** Typical convergence for the different parameterizations.

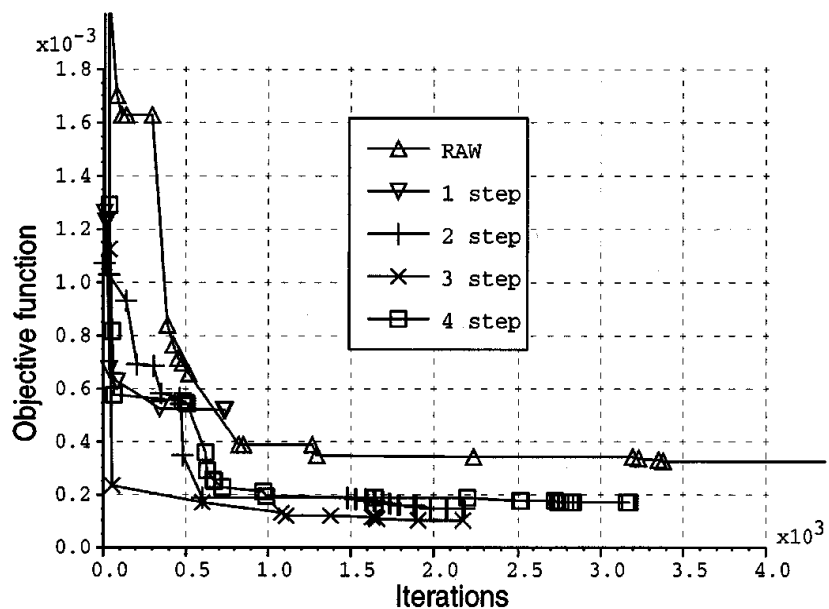

**Figure 6.** Convergence for another set of runs.

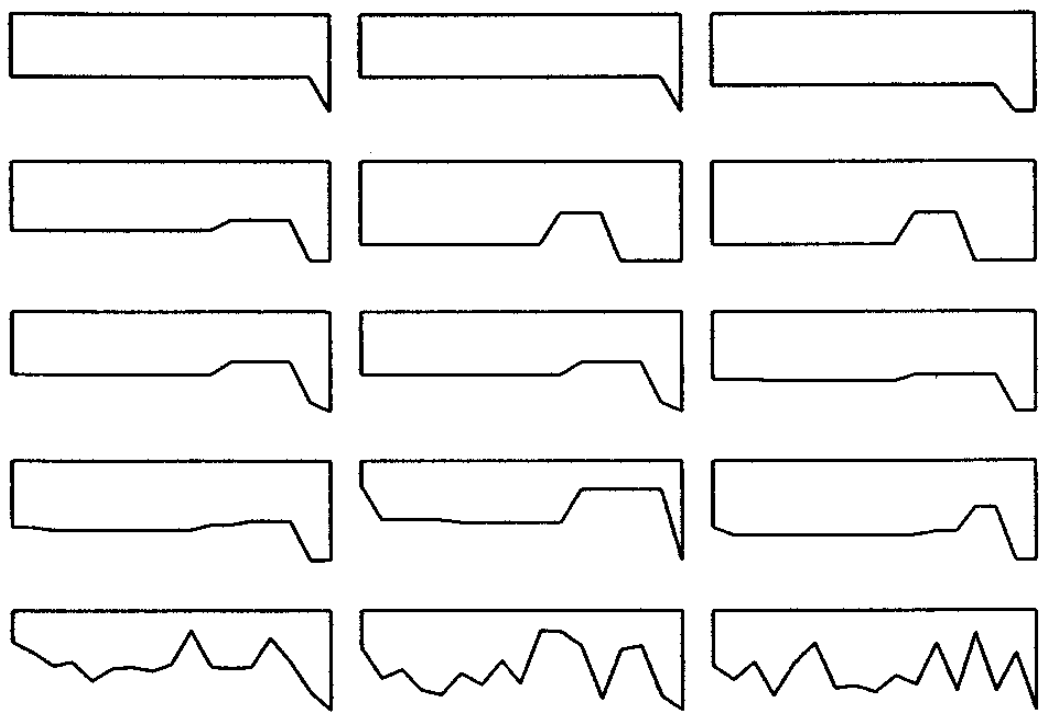

**Figure 7.** Typical pole shapes produced by the optimization process.

## IV. CONCLUSION

As expected the results clearly show that too many parameters increased the search time making an optimum shape difficult to find. Also (as expected) the 1 step parameterization was unable to achieve the highest field uniformity because the search space was too constrained. The 2 step parameterization performed well as did the 3 step parameterization. It can be seen from the results that the 3 pole parameterization was able to find a better configuration than the 2 step case but it was not so reliable.

**REFERENCES**

[1] R. J. Hill-Cottingham, J. F. Eastham, and I. R. Young, "Use of ironpowder material for MR-imaging magnet pole faces," *IEEE Transactions on Magnetics*, vol. 31, no. 6, pp. 4074–4076, November 1995.

[2] T. Miyamoto, H. Sakuria, H. Takabayashi, and M. Aoki, "A development of a permanent magnet assembly for MRI devices using Nd–Fe–B material," *IEEE Trans. Magn.*, vol. 25, no. 5, pp. 3907–3909, 1989.

[3] F. Glover, "Future paths for integer programming and links to artificial intelligence," *Computers and Operations Research*, vol. 5, pp. 533–549, 1986.

[4] R. Hooke and T. A. Jeeves, "Direct search solution of numerical problem," *Journal of Computing Machinery*, no. 8, pp. 212–229, 1961.

[5] M. de la Maza and D. Yuret, "Dynamic hill climbing," *AI Expert*, vol. 9, no. 3, pp. 26–31, 1994.